\documentclass[12pt]{article}
\usepackage[T2A]{fontenc}
\usepackage[cp1251]{inputenc}
\usepackage{amsfonts}
\usepackage{enumerate}

\begin{document}

\begin{center}
\large \textbf{Characterization of distributions of $Q$-independent\\
random variables on locally compact Abelian groups}

\bigskip

{\bf Margaryta Myronyuk}

\end{center}

 \makebox[20mm]{ }\parbox{125mm}{ \small Let  $X$ be a second countable locally compact
Abelian group.  We prove some group analogues
of the Skitovich--Darmois,   Heyde  and   Kac--Bernstein characterisation theorems  for
$Q$-independent random variables taking  values in the  group $X$. The proofs of
these theorems are reduced to solving some functional equations on
the character group of the group $X$.} 

\bigskip


\bigskip


\bigskip

\textbf{1. Introduction}

\bigskip

In the article  \cite{KS} A.M. Kagan and G.J. Sz\'ekely introduced a
notion of $Q$-independence of random variables which generalizes a
notion of independence of random variables. They proved that the Cram\'er
theorem about decomposition of Gaussian distribution  holds true if instead of independence
$Q$-independence is considered. They also proved that the
Skitovich--Darmois, 
Marcinkiewicz   and Vershik theorems hold true for $Q$-independent random
variables. In the papers \cite{Rao1} and \cite{Rao2} B.L.S. Prakasa Rao
considered some generalizations of the Kotlarski theorem and similar
results (see \cite{Kotlarski}) for $Q$-independent random variables.

The series of characterization theorems of
mathematical statistics were generalized on different algebraic
structures. Especially, much attention has been devoted to the study
of analogues of characterization theorems for locally compact
Abelian groups (see e.g. (\cite{Fe10}--\cite{FeBook2}). Following   A.M. Kagan and G.J. Sz\'ekely \cite{KS} in the paper
\cite{Fe2017} G.M. Feldman in a natural way introduced the notion of
$Q$-independence for random variables taking values in a locally
compact Abelian group. He proved that if we consider  $Q$-independence instead of independence, then the group analogue of the Cram\'er theorem (\cite{Fe1}, \cite{GF4}), and some group analogues of the the  Skitovich--Darmois (\cite{Fe16})  and  Heyde theorems (\cite{Fe30})
hold true for the same classes of groups. 

In this paper we continue research of characterization theorems for
$Q$-independent random variables with values in a locally compact Abelian group $X$ which were
started in \cite{Fe2017}. We study a group analog of the
Skitovich--Darmois theorem for $Q$-independent random variables and
linear forms with integer coefficients (\cite{Fe10}). We also study
a group analog of the Heyde theorem for $Q$-independent random
variables,  where coefficients of linear forms   are topological
automorphisms of the group $X$ (\cite{Fe18}). We prove that the
class of groups $X$ where  these theorems hold true does not change if
we replace independence for $Q$-independence. In contrast to these
results we show that the class of groups $X$ where the 
Kac--Bernstein theorem (\cite{Fe5}) holds true narrows if we replace
the condition of independence  for
$Q$-independence.

In the paper we suppose that $X$ is a second countable locally compact
Abelian group. Denote by $Y$ the character group of $X$, and by
$(x,y)$ the value of a character $y \in Y$ at $x\in X$. If $K$ is a closed subgroup of $X$, denote by
 $A(Y, K) = \{y \in Y: (x, y) = 1$ \mbox{ for all } $x \in K \}$
its annihilator. Denote by
${\rm Aut}(X)$ the group of topological automorphisms of the group
$X$, and by $I$ the identity automorphism of a group. If
$\alpha\in{\rm Aut}(X)$, then the adjoint  automorphism
$\widetilde\alpha\in{\rm Aut}(Y)$ is defined as follows $(x,
\widetilde\alpha y) = (\alpha x, y)$ for all $x \in X$, $y \in Y$.
Note that  $\alpha\in {\rm Aut}(X)$ if and only if
$\widetilde\alpha\in {\rm Aut}(Y)$. Denote by $\mathbb{R}$   the
group of real numbers, by $\mathbb{T}=\{z\in \mathbb{C}: |z|=1\}$
the circle group    (the one dimensional torus), and by $\mathbb{Z}$
the group of integers. Let $n$ be an integer. Denote by $f_n$ the
mapping of $X$ into $X$ defined by the formula $f_nx=nx$. Put
$X_{(n)}={\rm Ker} f_n$ and $X^{(n)}=f_n(X)$. A group $X$ is called
a Corwin group if $X^{(2)}=X$.

If $\xi$ is a  random variable with values in the group $X$, then
denote by $\mu_\xi$ its distribution and by
$$
\widehat\mu_\xi(y)={\bf E}[(\xi, y)]=\int_{X}(x, y)d \mu_\xi(x),
\quad y\in Y,$$ the characteristic function of the distribution
$\mu_\xi$.

Let $f(y)$ be a function on the group $Y$, and let $h \in Y$. Denote
by   $\Delta_h$   the finite difference operator
$$
\Delta_h f(y)=f(y+h)-f(y).
$$
We remind that a function $f(y)$ on $Y$ is called a polynomial  if
$$\Delta_{h}^{n+1}f(y)=0$$ for some   $n$ and for all $y, h \in Y$.
The minimal $n$ for which this equality holds is called   the degree
of the polynomial $f(y)$.

Let $\xi_1, \dots, \xi_n$ be random variables with values in the
group $X$. Following \cite{Fe2017} we say that the random variables
$\xi_1, \dots, \xi_n$ are $Q$-independent if the characteristic
function of the vector $(\xi_1, \dots, \xi_n)$ can be represented in
the form
\begin{equation}\label{i0}
    \widehat\mu_{(\xi_1, \dots, \xi_n)}(y_1, \dots, y_n)={\bf E}[(\xi_1,
y_1)\cdots(\xi_n,
y_n)]=$$$$=\left(\prod_{j=1}^n\widehat\mu_{\xi_j}(y_j)\right)\exp\{q(y_1,
\dots, y_n)\}, \quad y_j\in Y,
\end{equation}
where $q(y_1, \dots, y_n)$ is a continuous polynomial on the group
$Y^n$. We will also assume that $q(0, \dots, 0)=0$.


Denote by ${\rm M}^1(X)$ the convolution semigroup of probability
distributions on the group $X$. We remind that a distribution
$\gamma\in {\rm M}^1(X)$  is called Gaussian (see \cite[Chapter
IV]{Pa})  if its characteristic function is represented in the form
\begin{equation}\label{i00}
\widehat\gamma(y)=(x,y)\exp\{-\varphi(y)\}, \quad y\in Y,
\end{equation}
where $x \in X$, and $\varphi(y)$ is a continuous non-negative
function on the group $Y$ satisfying the equation
\begin{equation}\label{i1}
    \varphi(u+v)+\varphi(u-v)=2[\varphi(u)+\varphi(v)],
    \quad u,
    v\in Y.
\end{equation}
Denote by $\Gamma(X)$ the set of Gaussian distributions on the group $X$. We
note that according this definition the generated distributions are
Gaussian.  

\bigskip

\textbf{2. The Skitovich-Darmois theorem
for Q-independent random variables}

\bigskip

V. Skitovich and G. Darmois proved one of the most famous characterization
theorem of mathematical statistics in which the Gaussian
distribution on the real line is characterized by the independence of
two linear forms of $n$ independent random variables
(\cite[Ch.
3]{Kag-Lin-Rao}). A group analogue of the Skitovich--Darmois theorem
for random variables with values in a locally compact Abelian group, where coefficients of linear forms are integer,  was considered in \cite{Fe10} (see
also \cite[\S 10.7]{FeBook2}). We need the following definition. A
set of integers $\{a_j\}$ is said to be admissible for a group $X$
if $X^{(a_j)} \ne \{0\}$ for all $j$. Let $\xi_j, \ j = 1, 2,\dots,
n, \ n \ge 2,$ be independent random variables with values in $X$.
The admissibility of a set $\{a_j\}_{j=1}^n$ when considering the
linear form $L = a_1\xi_1 + \dots + a_n\xi_n$ is a group analogue of
the condition $a_j \ne 0$, $j = 1, 2,\dots, n$, for the case  of $X
= \mathbb{R}$. The following result holds.

\medskip

\textbf{Theorem A.} \textit{Let $X$ be a second countable locally
compact Abelian group. Let $\{a_j\}_{j=1}^n$ and $\{b_j\}_{j=1}^n$
be admissible sets of integers for the group $X$. Let
$\xi_1,\dots,\xi_n$ be independent random variables with values in
$X$ and distributions $\mu_j$ with non-vanishing characteristic
functions. The independence of the linear forms $L_1 = a_1\xi_1 +
\dots + a_n\xi_n$ and $L_2 = b_1\xi_1 + \dots + b_n\xi_n$ yields
that all $\mu_j \in \Gamma(X)$ if and only if either  $X$ is a
torsion-free group or $X^{(p)} = \{0\}$, where $p$ is a prime number
$($in the last case all $\mu_j$ are degenerate distributions$)$.}

\medskip

We prove that Theorem A remains true if we change the condition of
independence of random variables $\xi_1,\dots,\xi_n$ and linear
forms $L_1$ and $L_2$ for $Q$-independence. The following statement
is valid.

\medskip

\textbf{Theorem 1.} \textit{Let $X$ be a second countable locally
compact Abelian group. Let $\{a_j\}_{j=1}^n$ and $\{b_j\}_{j=1}^n$
be admissible sets of integers for the group $X$. Let
$\xi_1,\dots,\xi_n$ be Q-independent random variables with values in
$X$ and distributions $\mu_j$ with non-vanishing characteristic
functions. The Q-independence of the linear forms $L_1 = a_1\xi_1 +
\dots + a_n\xi_n$ and $L_2 = b_1\xi_1 + \dots + b_n\xi_n$ yields
that all $\mu_j \in \Gamma(X)$ if and only if either  $X$ is a
torsion-free group or  $X^{(p)} = \{0\}$, where $p$ is a prime
number  (in the last case all $\mu_j$ are degenerate
distributions).}

\medskip

To prove Theorem 1 we need some lemmas. The proof of the following
lemma almost literally coincides with the similar lemma of the paper
\cite{Fe2017}.

\medskip

\textbf{Lemma 1.} \textit{Let $X$ be a second countable locally
compact Abelian group and $Y$ be its character group. Let
$\{a_j\}_{j=1}^n$ and $\{b_j\}_{j=1}^n$ be sets of integers. Let
$\xi_1,\dots,\xi_n$ be Q-independent random variables with values in
$X$ and distributions $\mu_j$. The linear forms $L_1 = a_1\xi_1 +
\dots + a_n\xi_n$ and $L_2 = b_1\xi_1 + \dots + b_n\xi_n$ are
$Q$-independent if and only if the characteristic functions
$\widehat\mu_j(y)$ satisfy the equation}
\begin{equation}\label{l1}
    \prod_{i=1}^n\widehat\mu_j(a_j u+ b_j
    v) =\left(\prod_{i=1}^n\widehat\mu_j( a_j
    u)\prod_{i=1}^n\widehat\mu_j( b_j v)\right) \exp\{q(u,v)\}, \quad u, v \in
    Y,
\end{equation}
\textit{where $q(u,v)$ is a continuous polynomial on the group
$Y^2$, $q(0,0)=0$.}

\medskip

\textbf{Proof}. On the one hand, since the random variables $\xi_1, \dots, \xi_n$ are
are $Q$-independent, the characteristic function of the vector
$(L_1,L_2)$ is of the form
\begin{equation}\label{l1.1}
    \widehat\mu_{(L_1, L_2)}(u, v)={\bf E}[(L_1, u)(L_2, v)]={\bf
    E}[(a_1\xi_1+\cdots+a_n\xi_n, u)(b_1\xi_1+\cdots+ b_n\xi_n,
    v)]=$$$$={\bf E}[(\xi_1, a_1u+b_1v)\cdots(\xi_n, a_nu+b_nv)]=$$$$=
    \left(\prod_{j=1}^n\widehat\mu_{\xi_j}(a_ju+
    b_jv)\right)\exp\{q_1(a_1u+b_1v,\dots, a_nu+b_nv)\}, \quad u, v\in
    Y,
\end{equation}
where $q_1(y_1, \dots, y_n)$ is a continuous polynomial on the group
$Y^n$.

On the other hand, it follows from the $Q$-independence of $\xi_1,
\dots, \xi_n$ that
\begin{equation}\label{l1.2}
    \widehat\mu_{L_1}(y)={\bf E}[(L_1, y)]={\bf
    E}[(a_1\xi_1+\cdots+a_n\xi_n, y)]=
    \left(\prod_{j=1}^n\widehat\mu_{\xi_j}(a_jy)\right)
\exp\{q_1(a_1y,\dots, a_ny)\}, \quad y\in Y,
\end{equation}
\begin{equation}\label{l1.3}
    \widehat\mu_{L_2}(y)={\bf E}[(L_2, y)]={\bf
    E}[(b_1\xi_1+\cdots+b_n\xi_n, y)]=
    \left(\prod_{j=1}^n\widehat\mu_{\xi_j}(b_jy)\right)
    \exp\{q_1(b_1y,\dots, b_ny)\}, \quad y\in Y.
\end{equation}
Assume that the linear forms $L_1$ and $L_2$ are $Q$-independent.
Then the characteristic function of the vector $(L_1,L_2)$ can be
written in the form
\begin{equation}\label{l1.4}
    \widehat\mu_{(L_1, L_2)}(u, v)=\widehat\mu_{L_1}(u)\widehat\mu_{L_2}(v)
    \exp\{q_2(u, v)\}, \quad u, v\in Y,
\end{equation}
where $q_2(u, v)$ is a continuous polynomial on the group $Y^2$. Put
\begin{equation}\label{l1.5}
    q(u, v)=-q_1(a_1u+b_1v,\dots,a_nu+b_nv)+q_1(a_1u,\dots,
    a_nu)+q_1(b_1v,\dots, b_nv)+q_2(u,v).
\end{equation}
It follows from the definition of a polynomial on a group   that
$q(u, v)$ is a continuous polynomial on the group $Y^2$. Obviously,
(\ref{l1}) follows from (\ref{l1.1})--(\ref{l1.5}). If (\ref{l1})
holds, then (\ref{l1.4}) follows from (\ref{l1.1})--(\ref{l1.3}),
where $q_2(u, v)$ is defined by formula (\ref{l1.5}). Obviously,
$q_2(u, v)$ is a continuous polynomial on the group $Y^2$. Lemma 1
is proved. $\Box$

\medskip

\textbf{Lemma 2} (\cite{Fe1}, \cite{GF4}). \textit{Let $X$ be a
second countable locally compact Abelian group. Assume that $X$
contains no subgroup topologically isomorphic to the circle group
$\mathbb{T}$. Let $\xi_1$ and $\xi_2$ be independent random
variables with values in the group $X$. If the random variable
$\xi=\xi_1+\xi_2$ has a Gaussian distribution, then  $\xi_j$, $j=1,
2,$ are also Gaussian.}

\medskip

\textbf{Lemma 3} (\cite{Fe6}). \textit{Let $X$ be a second countable
locally compact Abelian group and $Y$ be its character group. Assume
that $X$ contains no subgroup topologically isomorphic to the circle
group $\mathbb{T}$. Let $f(y)$ be a characteristic function on the
group $Y$. If $f(y)$ is of the form
$$
f(y)=\exp\{P(y)\}, \quad  y\in Y,
$$
where $P(y)$ is a continuous polynomial, then $P(y)$ is a polynomial
of degree $\le 2$, and $f(y)$ is the characteristic function of a
Gaussian distribution on the group $X$.}

\medskip

\textbf{Proof of Theorem 1}. \textbf{Sufficiency.} By Lemma 1, it follows from the conditions of the theorem
that the characteristic functions $\widehat\mu_j(y)$ satisfy
equation (\ref{l1}). Put $\nu_j = \mu_j
* \bar \mu_j$. Then $\widehat \nu_j(y) = |\widehat \mu_j(y)|^2 > 0,$
\ $y \in Y$. Obviously, the characteristic functions $\widehat
\nu_j(y)$
 satisfy the equation
\begin{equation}\label{t1}
    \prod_{i=1}^n\widehat\nu_j(a_j u+ b_j
    v) =\left(\prod_{i=1}^n\widehat\nu_j( a_j
    u)\prod_{i=1}^n\widehat\nu_j( b_j v)\right) \exp\{r(u,v)\}, \quad u, v \in
    Y,
\end{equation}
where $r(u,v)=q(u,v)+\overline{q(u, v)}$ is a real valued continuous polynomial on the
group $Y^2$, $r(0,0)=0$.

Assume first that $X$ is a torsion-free group. Then the group $X$
contains no subgroup topologically isomorphic to the circle group
$\mathbb{T}$. Therefore, if we prove that $\nu_j \in \Gamma(X)$,
then applying Lemma 2 we obtain that $\mu_j \in \Gamma(X)$. Hence we
can assume from the beginning that $\widehat \mu_j(y)
> 0$, $y \in Y$, $j = 1, 2, \dots, n$, $n \ge 2$.

Set $\varphi_j(y)=-\ln \widehat\mu_j(y)$. We conclude from
(\ref{l1}) that
\begin{equation}\label{t2}
    \sum_{j = 1}^{n}{\varphi_j(a_j u + b_j v)} = P(u)+Q(v)+q(u,v),
\quad u,v\in Y,
\end{equation}
where
\begin{equation}\label{t3}
    P(u)=\sum_{j = 1}^{n}{\varphi_j(a_j u)},\quad Q(v)=\sum_{j =
1}^{n}{\varphi_j(b_j v)}.
\end{equation}

We use the finite difference method to solve equation (\ref{t2}).
Let $h_n$ be an arbitrary element of the group $Y$. Substitute
$u+b_n h_n$ for $u$ and $v-a_n h_n$ for $v$  in equation (\ref{t2}).
Subtracting equation (\ref{t2}) from the resulting equation we
obtain
\begin{equation}\label{t4}
    \sum_{j = 1}^{n-1} \Delta_{l_{n,j}}{\varphi_j(a_j u + b_j v)}
    =\Delta_{b_n h_{n}} P(u)+\Delta_{-a_n h_{n}} Q(v) +\Delta_{(b_n h_n,-a_n h_{n})} q(u,v),
\quad u,v\in Y,
\end{equation}
where $l_{n,j}=(a_j b_n -b_j a_n)h_n$, $j=1, 2,\dots,n-1$. Let
$h_{n-1}$ be an arbitrary element of the group $Y$. Substitute
$u+b_{n-1}h_{n-1}$ for $u$ and $v-a_{n-1}h_{n-1}$ for $v$  in
equation (\ref{t4}). Subtracting equation (\ref{t4}) from the
resulting equation we obtain
\begin{equation}\label{t5}
    \sum_{j = 1}^{n-2} \Delta_{l_{n-1,j}} \Delta_{l_{n,j}}{\varphi_j(a_j u + b_j v)}
    =\Delta_{b_{n-1}h_{n-1}} \Delta_{b_n h_{n}} P(u)+
    \Delta_{-a_{n-1}h_{n-1}} \Delta_{-a_n h_{n}} Q(v)+$$
    $$+\Delta_{(b_{n-1}h_{n-1},-a_{n-1}h_{n-1})}\Delta_{(b_n h_n,-a_n h_{n})} q(u,v),
    \quad u,v\in Y,
\end{equation}
where $l_{n-1,j}= (a_j b_{n-1}-b_j a_{n-1})h_{n-1}$, $j=1,
2,\dots,n-2$. Arguing as above we get the equation
\begin{equation}\label{t6}
    \Delta_{l_{2,1}} \Delta_{l_{3,1}}\dots \Delta_{l_{n,1}}{\varphi_1(a_1 u + b_1 v)}=$$
    $$=\Delta_{b_2 h_2}\Delta_{b_3 h_3}\dots
    \Delta_{b_n h_{n}} P(u)+\Delta_{-a_{2}h_{2}} \Delta_{-a_{3}h_{3}}\dots
    \Delta_{-a_n h_{n}} Q(v)+$$
    $$+\Delta_{(b_2 h_2,-a_{2}h_{2})} \Delta_{(b_3 h_3,-a_{3}h_{3})}\dots\Delta_{(b_n h_n,-a_n h_{n})} q(u,v),
\quad u,v\in Y,
\end{equation}
where $h_m$ is an arbitrary element of the group $Y$, $l_{m,1}= (a_1
b_m-b_1 a_m)h_m$, $m=2, 3, \dots,n$.

Let $h_{1}$ be an arbitrary element of the group $Y$. Substitute
$u+b_{1}h_{1}$ for $u$ and $v-a_{1}h_{1}$ for $v$ in equation
(\ref{t6}). Subtracting equation (\ref{t6}) from the resulting
equation we obtain
\begin{equation}\label{t7}
   \Delta_{b_1 h_1}\Delta_{b_2 h_2}\Delta_{b_3 h_3}\dots
    \Delta_{b_n h_{n}} P(u)+\Delta_{-a_{1}h_{1}}
    \Delta_{-a_{2}h_{2}} \Delta_{-a_{3}h_{3}}\dots \Delta_{-a_n h_{n}} Q(v)+$$
    $$+\Delta_{(b_1 h_1,-a_{1}h_{1})} \Delta_{(b_2 h_2,-a_{2}h_{2})} \Delta_{(b_3 h_3,-a_{3}h_{3})}\dots\Delta_{(b_n h_n,-a_n h_{n})} q(u,v)=0,
\quad u,v\in Y.
\end{equation}
Let $h$ be an arbitrary element of the group $Y$. Substitute $u+h$
for $u$  in equation (\ref{t7}) and subtract equation (\ref{t7})
from the resulting equation. We get
\begin{equation}\label{t8}
   \Delta_{h}\Delta_{b_1 h_1}\Delta_{b_2 h_2}\Delta_{b_3 h_3}\dots
    \Delta_{b_n h_{n}} P(u)+$$
    $$+\Delta_{(h,0)}\Delta_{(b_1 h_1,-a_{1}h_{1})} \Delta_{(b_2 h_2,-a_{2}h_{2})} \Delta_{(b_3 h_3,-a_{3}h_{3})}\dots\Delta_{(b_n h_n,-a_n h_{n})} q(u,v)=0,
\quad u\in Y.
\end{equation}

Note that if $h$ and $k$ are arbitrary elements of the group
$Y$, by the condition
\begin{equation}\label{t9}
  \Delta^{l+1}_{(h, k)}q(u, v)=0,
\quad u, v \in Y,
\end{equation}
for some $l$.

We assumed that  $X$ is a torsion-free group, i.e. $X_{(m)}=\{0\}$
for all $m\in \mathbb{Z}$, $m\ne 0$. This implies that
$\overline{Y^{(m)}}=Y$ for all $m\in \mathbb{Z}$, $m\ne 0$. In
particular, $\overline{Y^{(a_j)}}=Y$, $\overline{Y^{(b_j)}}=Y, \ j=1,
2, \dots, n$. Hence taking into account that $h$ and $h_j$ are arbitrary
elements of the group $Y$, it follows from (\ref{t8}) that
\begin{equation}\label{t10}
   \Delta_{h}^{n+1}P(u)
   +\Delta_{(h,0)}\Delta_{(h,k)}^{n} q(u,v)=0,
\quad u,v\in Y.
\end{equation}

We apply to (\ref{t10}) the operator $\Delta^{l+1}_{(h, k)}$.
Taking into account (\ref{t9}), we get
$$
   \Delta_{h}^{l+n+2}P(u)=0, \quad u\in Y.
$$

Thus, $P(y)$ is a continuous polynomial on the group $Y$. Consider the
distribution ${\gamma=f_{a_1}(\mu_1)*\dots* f_{a_n}(\mu_n)}$. Then
$\widehat{f_{a_j}(\mu_j)}(y)=\widehat\mu(a_jy), \ {j=1, 2,\dots,n}$.
We get that
$$
    \widehat\gamma(y)=\prod_{j = 1}^{n}\widehat\mu_j(a_j y), \quad y\in Y.
$$
Taking into account (\ref{t3}), we have
$$\widehat\gamma(y)=e^{-P(y)}.$$ Since $X$ is a torsion free group, the group $X$
contains no subgroup topologically isomorphic to the circle group
$\mathbb{T}$. Hence Lemma 3 implies that $\gamma\in\Gamma(X)$. Then
by Lemma 2 we obtain that all $f_{a_j}(\mu_j)\in\Gamma(X)$. This
implies that every function
 $\varphi_j(y)$ satisfies  equation (\ref{i1}) on the set
${Y^{(a_j)}}$. Since $\overline{Y^{(a_j)}}=Y$, $j=1, 2\dots, n$,
 every function $\varphi_j(y)$ satisfies  equation
(\ref{i1}) on the group $Y$. This means that all
$\mu_j\in\Gamma(X)$.

Let $X^{(p)} = \{0\}$, where $p$ is a prime number. Then it follows
from admissibility of the sets  $\{a_j\}_{j=1}^n$ and
$\{b_j\}_{j=1}^n$ for the group $X$ that $f_{a_j}, f_{b_j}\in{\rm
Aut}(X)$, $j=1, 2,\dots, n$. We can apply Theorem 1 of \cite{Fe2017} and get that all distributions $\mu_j$ are Gaussian. Note that the condition $X^{(p)} = \{0\}$ implies that the connected component of zero of the group $X$ is equal to $\{0\}$. Taking into account that the support of a Gaussian distribution is a coset of the connected component of zero of the group $X$, it means that
  all $\mu_j$ are degenerate
distributions.

\textbf{Necessity.} Since Q-independent random variables are
independent, necessity follows from Theorem A. $\Box$

\bigskip

\textbf{3. The Heyde theorem for Q-independent random variables}

\bigskip

The similar result to the Skitovich--Darmois theorem  was proved by C.C. Heyde
  (\cite[\S 13.4]{Kag-Lin-Rao}), where the condition
of the independence of two linear forms is replaced by the condition
of the symmetry of the conditional distribution of the linear form
$L_2$ given $L_1$. The following group analogue of the
Heyde theorem was proved in \cite{Fe18} (see also \cite[\S
16.2]{FeBook2}).

\medskip

\textbf{Theorem B.} \textit{Let $X$ be a second countable locally
compact Abelian group. Let $\alpha_j, \ \beta_j, \ j = 1, 2,\dots,
n, \ n \ge 2,$ be topological automorphisms of $X$ such that
$\beta_i\alpha_i^{-1} \pm \beta_j\alpha_j^{-1} \in {\rm Aut}(X)$ \
for all $i \ne j$. Let $\xi_1,\dots,\xi_n$ be independent random
variables with values in $X$ and distributions $\mu_j$ with
non-vanishing characteristic functions. The symmetry of the
conditional distribution of the linear form $L_2 = \beta_1\xi_1 +
\dots + \beta_n\xi_n$ given $L_1 = \alpha_1\xi_1 + \dots +
\alpha_n\xi_n$ implies that all $\mu_j \in \Gamma(X)$ if and only if
$X$ contains no elements of order $2$.}

\medskip

We prove that Theorem B remains true if we change the condition of
independence of random variables $\xi_1,\dots,\xi_n$ for
$Q$-independence. The following statement is valid.

\medskip

\textbf{Theorem 2.} \textit{Let $X$ be a second countable locally
compact Abelian group. Let $\alpha_j, \ \beta_j, \ j = 1, 2,\dots,
n, \ n \ge 2,$ be topological automorphisms of $X$ such that
$\beta_i\alpha_i^{-1} \pm \beta_j\alpha_j^{-1} \in {\rm Aut}(X)$ \
for all $i \ne j$. Let $\xi_1,\dots,\xi_n$ be Q-independent random
variables with values in $X$ and distributions $\mu_j$ with
non-vanishing characteristic functions. The symmetry of the
conditional distribution of the linear form $L_2 = \beta_1\xi_1 +
\dots + \beta_n\xi_n$ given $L_1 = \alpha_1\xi_1 + \dots +
\alpha_n\xi_n$ implies that all $\mu_j \in \Gamma(X)$ if and only if
$X$ contains no elements of order $2$.}

\medskip

To prove Theorem 2 we need the following lemma.

\medskip

\textbf{Lemma 4.} \textit{Let  $X$ be a second countable locally
compact Abelian group, $Y$ be its character group. Let $\alpha_j, \
\beta_j, \ j = 1, 2,\dots, n, \ n \ge 2,$ be topological
automorphisms of the group $X$. Let  $\xi_1,\dots,\xi_n$ be
$Q$-independent random variables with values in $X$ and
distributions $\mu_j$. The conditional distribution of the linear
form $L_2 = \beta_1\xi_1 + \dots + \beta_n\xi_n$ given $L_1 =
\alpha_1\xi_1 + \dots + \alpha_n\xi_n$ is symmetric if and only if
the characteristic functions $\widehat\mu_{\xi_j}(y)$  satisfy the
equation
\begin{equation}\label{l4}
    \prod_{j=1}^n \widehat\mu_j(\widetilde\alpha_j u+\widetilde\beta_j v )=
    \prod_{j=1}^n \widehat\mu_j(\widetilde\alpha_j u-\widetilde\beta_j
    v) \exp\{q(u,v)\}, \quad u, v \in Y,
\end{equation}
where $q(u, v)$ is a continuous polynomial on the group $Y^2$, $q(0,
0)=0$.}

\medskip

\textbf{ Proof}. Since the random variables $\xi_1,\dots,\xi_n$ are
$Q$-independent,   the characteristic function of the vector
$(L_1,L_2)$ is of the form
\begin{equation}\label{l4.1}
    \widehat\mu_{(L_1, L_2)}(u, v)={\bf E}[(L_1, u)(L_2, v)]=
    {\bf E}[(\alpha_1\xi_1 + \dots + \alpha_n\xi_n, u)(\beta_1\xi_1 + \dots + \beta_n\xi_n, v)]$$
    $$={\bf E}[(\xi_1, \widetilde{\alpha_1}u+ \widetilde{\beta_1}v)\dots
    (\xi_n, \widetilde{\alpha_n}u+ \widetilde{\beta_n}v)]$$$$= \prod_{j=1}^n \widehat\mu_j(\widetilde\alpha_j u+\widetilde\beta_j v )
    \exp\{q_1(\widetilde{\alpha_1}u+ \widetilde{\beta_1}v,\dots, \widetilde{\alpha_n}u+ \widetilde{\beta_n}v)\}, \quad u, v\in Y,
\end{equation}
where $q_1(y_1, \dots, y_n)$ is a continuous polynomial on the group
$Y^n$, $q_1(0,\dots, 0)=0$. Similarly, the characteristic function
of the vector $(L_1,-L_2)$ is of the form
\begin{equation}\label{l4.2}
    \widehat\mu_{(L_1, -L_2)}(u, v)= \prod_{j=1}^n \widehat\mu_j(\widetilde\alpha_j u-\widetilde\beta_j v )
    \exp\{q_1(\widetilde{\alpha_1}u- \widetilde{\beta_1}v,\dots, \widetilde{\alpha_n}u- \widetilde{\beta_n}v)\}, \quad u, v\in
    Y.
\end{equation}
The symmetry of the conditional distribution of the linear form
$L_2$ given  $L_1$ means that the random vectors $(L_1, L_2)$ and
$(L_1, -L_2)$ are identically distributed, i.e. the  
characteristic functions (\ref{l4.1}) and (\ref{l4.2}) are the same.
Put $q(u, v)=-q_1(\widetilde{\alpha_1}u+ \widetilde{\beta_1}v,\dots,
\widetilde{\alpha_n}u+
\widetilde{\beta_n}v)+q_1(\widetilde{\alpha_1}u-
\widetilde{\beta_1}v,\dots, \widetilde{\alpha_n}u-
\widetilde{\beta_n}v)$. It is obvious that $q(u, v)$ is a continuous
polynomial on the group $Y^2$, and (\ref{l4}) follows from
(\ref{l4.1}) and (\ref{l4.2}). $\Box$

\medskip

\textbf{Proof of Theorem 2}. \textbf{Sufficiency.} Let $Y$ be the
character group of the group $X$. We can put $\zeta_j = \alpha_j \xi_j$ and
reduce the proof of the theorem to the case when $L_1 = \xi_1 +
\dots + \xi_n$ and $L_2 =\delta_1 \xi_1 + \dots + \delta_n\xi_n$,  
$\delta_j \in {\rm Aut}(X)$. The conditions $\beta_i\alpha_i^{-1} \pm
\beta_j\alpha_j^{-1} \in {\rm Aut}(X)$ for all $i \ne j$ are
transformed into the conditions $\delta_i \pm \delta_j \in {\rm
Aut}(X)$ for all $i \ne j$. By Lemma 4, the symmetry of the
conditional distribution of the linear form $L_2 =\delta_1 \xi_1 + \dots +
\delta_n\xi_n$ given $L_1 = \xi_1 + \dots + \xi_n$ implies that the
characteristic functions  $\widehat\mu_j(y)$ satisfy equation
(\ref{l4}) which takes the form
\begin{equation}\label{t2.2}
    \prod_{j=1}^n \widehat\mu_j(u+\widetilde\delta_j v )=
    \prod_{j=1}^n \widehat\mu_j(u-\widetilde\delta_j
    v) \exp\{q(u,v)\}, \quad u, v \in Y.
\end{equation}
Set $\nu_j = \mu_j * \bar \mu_j$. Then $\widehat \nu_j(y) =
|\widehat \mu_j(y)|^2 > 0, \ y \in Y$. Obviously, the characteristic
functions $\widehat \nu_j(y)$ satisfy equation
\begin{equation}\label{t2.1}
    \prod_{j=1}^n \widehat\nu_j(u+\widetilde\delta_j v )=
    \prod_{j=1}^n \widehat\nu_j(u-\widetilde\delta_j
    v) \exp\{r(u,v)\}, \quad u, v \in Y,
\end{equation}
where $r(u,v)=q(u,v)+\overline{q(u, v)}$ is a a real valued continuous polynomial on the
group $Y^2$, $r(0,0)=0$.

Set $\varphi_j(y) = -\ln \widehat\nu_j(y)$. It follows from
(\ref{t2.1}) that the functions $\varphi_j(y)$ satisfy equation
\begin{equation}\label{t2.3}
    \sum_{j = 1}^{n} [\varphi_j(u + \widetilde \delta_j
    v ) - \varphi_j(u - \widetilde \delta_j v )]=r(u,v), \quad u, v \in Y.
\end{equation}
We use the finite difference method to solve equation (\ref{t2.3}).

Let $k_1$ be an arbitrary element of the group $Y$. Set $h_1 =
\widetilde\delta_n k_1 $, then $h_1 -\widetilde\delta_n k_1  = 0$.
Substitute $u+h_1$ for $u$ and $v+k_1$ for $v$ in equation
(\ref{t2.3}). Subtracting equation (\ref{t2.3}) from the resulting
equation we obtain
\begin{equation}\label{t2.4}
    \sum_{j = 1}^{n} \Delta_{l_{1,j}}\varphi_j(u +
    \widetilde \delta_j v ) - \sum_{j =1}^{n-1}
    \Delta_{l_{1,j+n}}\varphi_j(u - \widetilde \delta_j v)= \Delta_{(h_1,k_1)} r(u,v), \quad u, v
    \in Y,
\end{equation}
where $l_{1,j} = h_1 + \widetilde\delta_jk_1 = (\widetilde\delta_n +
\widetilde\delta_j) k_1 , \ j = 1, 2,\dots, n, \ l_{1,j+n} = h_1 -
\widetilde\delta_j k_1  = (\widetilde\delta_n -\widetilde\delta_j)
k_1, \ j = 1, 2,\dots, n-1.$ Let $k_2$ be an arbitrary element of
the group $Y$. Set $h_2 = \widetilde\delta_{n-1} k_2 $. Then $h_2
-\widetilde\delta_{n-1} k_2  = 0$. Substitute $u+h_2$ for $u$ and
$v+k_2$ for $v$ in equation (\ref{t2.4}). Subtracting equation
(\ref{t2.4}) from the resulting equation we obtain
\begin{equation}\label{t2.5}
    \sum_{j = 1}^{n}
    \Delta_{l_{2,j}}\Delta_{l_{1,j}}\varphi_j(u + \widetilde \delta_j v)
    - \sum_{j =1}^{n-2} \Delta_{l_{2,j+n}}\Delta_{l_{1,j+n}}\varphi_j(u
    - \widetilde \delta_j v)=\Delta_{(h_2,k_2)}\Delta_{(h_1,k_1)} r(u,v), \quad u, v \in Y,
\end{equation}
where $l_{2,j} = h_2 + \widetilde\delta_j k_2  = (\widetilde\delta_n
+ \widetilde\delta_j) k_2, \ j = 1, 2,\dots, n, \ l_{2,j+n} = h_2 -
\widetilde\delta_j k_2 = (\widetilde\delta_n -\widetilde\delta_j)
k_2, \ j = 1, 2,\dots, n-2.$ Arguing as above in $n$ steps we get
the equation
\begin{equation}\label{t2.6}
    \sum_{j =
    1}^{n}\Delta_{l_{n,j}}\Delta_{l_{n-1,j}}\dots\Delta_{l_{1,j}}\varphi_j(u
    + \widetilde\delta_j v) = \Delta_{(h_n,k_n)}\dots\Delta_{(h_1,k_1)} r(u,v), \quad u, v \in Y,
\end{equation}
where $l_{p,j} = (\widetilde\delta_{n-p+1}+\widetilde\delta_j) k_p ,
\ p= 1, 2,\dots, n, \ j= 1, 2,\dots,n.$ Let $k_{n+1}$ be an
arbitrary element of the group $Y$. Set $h_{n+1} =
-\widetilde\delta_n k_{n+1}$, hence $h_{n+1} +\widetilde\delta_n
k_{n+1} = 0$. Substitute $u+h_{n+1}$ for $u$ and $v+k_{n+1}$ for $v$
in equation (\ref{t2.6}). Subtracting equation (\ref{t2.6}) from the
resulting equation we obtain
\begin{equation}\label{t2.7}
    \sum_{j =
    1}^{n-1}\Delta_{l_{n+1,j}}\Delta_{l_{n,j}}\Delta_{l_{n-1,j}}
    \dots\Delta_{l_{1,j}}\varphi_j(u + \widetilde\delta_j v ) = \Delta_{(h_{n+1},k_{n+1})}\dots\Delta_{(h_1,k_1)} r(u,v), \quad
    u, v \in Y,
\end{equation}
where $l_{n+1,j} = h_{n+1} + \widetilde\delta_j k_{n+1}  =
(\widetilde\delta_j - \widetilde\delta_n) k_{n+1}, \ j = 1, 2,\dots,
n-1.$ Equation (\ref{t2.7}) does not contain the function
$\varphi_n$. Arguing similarly we sequentially exclude the functions
$\varphi_{n-1}, \ \varphi_{n-2},\dots, \varphi_2$ from equation
(\ref{t2.7}). Finally we obtain
\begin{equation}\label{t2.8}
    \Delta_{l_{2n-1,1}}\Delta_{l_{2n-2,1}}\dots\Delta_{l_{1,1}}\varphi_1(u
    + \widetilde\delta_1 v) = \Delta_{(h_{2n-1},k_{2n-1})}\dots\Delta_{(h_1,k_1)} r(u,v), \quad u, v \in Y,
\end{equation}
where $l_{p,1} = (\widetilde\delta_1 +\widetilde\delta_{n-p+1}) k_p,
\ p= 1, 2, \dots, n-1, \ l_{n,1} = 2\widetilde\delta_1 k_n, \ l_{n +
p, 1} = (\widetilde\delta_1 - \widetilde\delta_{n-p+1}) k_{n + p}, \
p = 1, 2, \dots, n-1.$

Note that if $h$ and $k$ are arbitrary elements of $Y$, then
$$\Delta^{l+1}_{(h, k)}r(u, v)=0, \ u, v \in Y,$$  for some $l$.
Taking it into account and applying to the both sides of equation
(\ref{t2.8}) the operator $\Delta^{l+1}_{(h, k)}$, we obtain
\begin{equation}\label{t2.9}
    \Delta^{l+1}_{(h, k)}\Delta_{l_{2n-1,1}}\Delta_{l_{2n-2,1}}\dots\Delta_{l_{1,1}}\varphi_1(u
    + \widetilde\delta_1 v) = 0, \quad u, v \in Y.
\end{equation}

Taking into account that $k_p$ are arbitrary elements of the group
$Y$ and $\widetilde\delta_i \pm \widetilde\delta_j \in {\rm Aut}(Y)$
for $i \neq j$, we can substitute in (\ref{t2.9}) $l_{n,1} = 2k, \
l_{p,1} = h, \ p = 1, 2,\dots, n-1, n+1,\dots,2n-1$, where $k$ and
$h$ are arbitrary elements of the group $Y$. Substituting $v = 0$
into the resulting equation we obtain
\begin{equation}\label{t2.10}
    \Delta_{2k}\Delta_h^{2n+l - 1} \varphi_1(u) = 0, \quad u \in
    Y,
\end{equation}
where $k, h$ and $u$ are arbitrary elements of $Y$. Since the group
$X$ contains no elements of order 2, we have
$\overline{Y^{(2)}}=Y$. We deduce from (\ref{t2.10}) that
$$ \Delta_h^{2n+l} \varphi_1(u) = 0,
\quad u \in Y, $$ i.e. the function $\varphi_1(y)$ is a continuous
polynomial. Since the group $X$ contains no elements of order 2,
$X$ contains no subgroup topologically isomorphic to the circle
group ${\mathbb T}$. By Lemma 3, $\nu_1 \in \Gamma(X)$. Applying
Lemma 2 we get $\mu_1 \in \Gamma(X)$. Arguing as above we prove that
all  $\mu_j \in \Gamma(X)$.

\textbf{Necessity.} Since Q-independent random variables are
independent, necessity follows from Theorem B. $\Box$



\bigskip

\textbf{4. The Kac-Bernstein theorem for Q-independent identically
distributed random variables}

\bigskip

M. Kac   and S.N. Bernstein   proved that
the Gaussian distribution on the real line is characterized by the independence of the
sum and the difference of two independent random variables. The
following group analogue of the Kac--Bernstein theorem was proved
in \cite{Fe5} (see also \cite[\S 9.9]{FeBook2}).

\medskip

\textbf{Theorem C.} \textit{Let $X$ be a second countable locally
compact Abelian group. Let $\xi_1$ and $\xi_2$ be independent
identically distributed random variables with values in $X$ and
distribution $\mu$. The independence of the sum
$\xi_1+\xi_2$ and the difference $\xi_1-\xi_2$ yields that all $\mu\in
\Gamma(X)*I(X)$ if and only if the connected component of zero of
the group $X$ contains no more than one element of order $2$.}

\medskip

We prove that, unlike Theorems A and B, Theorem C does not remain
true if we change the condition of independence of random variables
$\xi_1 $ and $\xi_2 $ and the sum $\xi_1 + \xi_2$ and the difference
$\xi_1-\xi_2$ for Q-independence. In this case the distribution
$\mu\in\Gamma(X)*I(X)$ if and only if the connected component of
zero of the group $X$ contains no elements of order $2$. The following
statement is valid.

\medskip

\textbf{Theorem 3.} \textit{Let $X$ be a second countable locally
compact Abelian group. Let $\xi_1$ and $\xi_2$ be Q-independent
identically distributed random variables with values in $X$ and
distribution $\mu$. The Q-independence of the   sum
$\xi_1+\xi_2$ and the difference  $\xi_1-\xi_2$ yields that $\mu\in \Gamma(X)*I(X)$
if and only if the connected component of zero of the group $X$ contains no
elements of order $2$.}

\medskip

To prove Theorem 3 we need the following lemmas.

\medskip

\textbf{Lemma 5} (\cite{Fe2017}). \textit{Let $X$ be a second
countable locally compact Abelian group. Let $\xi_1$ and $\xi_2$ be
Q-independent random variables with values in $X$ and distributions
$\mu_j$. The Q-independence of the sum
$\xi_1+\xi_2$ and the difference 
$\xi_1-\xi_2$ yields that all $\mu_j\in \Gamma(X)*I(X)$ if and only
if the connected component of zero of the group $X$ contains no elements of
order $2$.}

\medskip

The following lemma was not formulated as a separate statement in
the paper \cite{Fe2017}, but it was essentially proven in it.

\medskip

\textbf{Lemma 6} (\cite{Fe2017}). \textit{Let $X=\mathbb{T}$. There
exist Q-independent identically distributed random variables $\xi_1$
and $\xi_2$ with values in $X$ and distribution $\mu$ with a
non-vanishing characteristic function such that the sum $\xi_1+\xi_2$ and the difference
$\xi_1-\xi_2$ are Q-independent, whereas $\mu\not\in \Gamma(X)$.}

\medskip

\textbf{Proof of Theorem 3.} \textbf{Sufficiency} follows from Lemma
5.

\textbf{Necessity.} Assume that the connected component of zero of
the group $X$ contains elements of order 2. Applying \cite[Lemma
7.6]{FeBook2} for $n=2$ we get that there exists a compact Corwin
subgroup $G$ of the group $X$ such that its character group is not a Corwin group.
Then by \cite[Lemma 7.7]{FeBook2} there exists a compact subgroup
$K$ of the group $X$ such that the factor-group $X/K$ contains a
subgroup $F$ topologically isomorphic to the circle group
$\mathbb{T}$. For this reason the distribution $\mu$ on the circle
group $\mathbb{T}$ from Lemma 6 can be considered as
a distribution  on the factor-group $X/K$. We will retain for it the
notation $\mu$. Since the character group of the factor-group $X/K$  is topologically isomorphic to the annihilator $A(Y, K)$, we can assume that the
characteristic function $\widehat{\mu}(y)$ is defined on $A(Y, K)$.
Consider on the group $Y$ the function
 $$
    h(y)=
    \left\{%
\begin{array}{ll}
    \widehat{\mu}(y), & \hbox{$y\in A(Y, K)$;} \\
    0, & \hbox{$y\notin A(Y, K)$.} \\
\end{array}%
\right.
 $$
Since $A(Y, K)$  is a subgroup and $\widehat{\mu}(y)$ is a positive
definite function, the function $h(y)$ is also positive definite
function (\cite[\S 2.12]{FeBook2}). Since $K$ is a compact group,
the annihilator $A(Y, K)$ is an open subgroup. Hence the function
$h(y)$ is continuous. By the Bochner theorem there exists
a distribution $\lambda\in{\rm M}^1(X)$ such that
$\widehat\lambda(y)=h(y)$. Let $\xi_j$ be independent identically
distributed random variables with values in the group $X$ and
distribution $\lambda$.

We verify that the sum $\xi_1+\xi_2$ and the difference $\xi_1-\xi_2$ are Q-independent.
Indeed, by Lemma 1, it suffices to show that the characteristic
function $\widehat\lambda(y)$ satisfy equation (\ref{l1}) which
takes the form
\begin{equation}\label{t3.1}
    \widehat\lambda(u+v) \widehat\lambda(u-v)=\widehat\lambda^2(u)\widehat\lambda(v)\widehat\lambda(-v) \exp\{q(u,v)\}, \quad u, v \in
    Y,
\end{equation}
where $q(u,v)$ is continuous polynomial on $Y^2$, $q(0,0)=0$.

Let $u, v \in A(Y, K)$. Then it is obvious that (\ref{t3.1}) holds,
because the function $\widehat{\mu}(y)$ satisfy equation
(\ref{t3.1}). If either $u \in  A(Y, K),  v \notin A(Y, K)$ or $v
\in A(Y, K), u \notin A(Y, K)$, then the both sides of equation
(\ref{t3.1}) are equal to zero. If $u, v \notin A(Y, K)$, then the
right-hand side of equation (\ref{t3.1}) is equal to zero. If the
left-hand side of equation (\ref{t3.1}) is not equal to zero, then
$u\pm v \in A(Y, K)$. This implies that $2u\in A(Y, K)$. Since $K$
is a compact Corwin group, applying \cite[Lemma 7.2]{FeBook2} we get
that $u\in A(Y, K)$. This contradicts the assumption. Thus the
left-hand side of (\ref{t3.1}) is also equal to zero. We proved that
the characteristic function $\widehat \lambda(y)$ satisfy equation
(\ref{t3.1}). Since $\mu\notin \Gamma(X/K)$, it is obvious that
$\lambda\notin \Gamma(X)*I(X)$. Necessity of Theorem 3 is proved.
$\Box$

\end{document}